\documentclass[10pt,a4paper]{article}
\usepackage[latin1]{inputenc}
\usepackage[english]{babel}
\usepackage{latexsym}
\usepackage{amsfonts}\usepackage{amsmath}

\newtheorem{theorem}{Theorem}
\newtheorem{proposition}[theorem]{Proposition}
\newtheorem{lemma}[theorem]{Lemma}
\newtheorem{corollary}[theorem]{Corollary}
\newtheorem{definition}[theorem]{Definition}
\newtheorem{remark}[theorem]{Remark}

\newtheorem{example}{Example}

\newenvironment{proof}{\noindent{\it Proof.\ }}{\hfill $\square$\vspace{2mm}}

\sectionmark \leftmark

\newcommand{\R}{\mathbb{R}}
\newcommand{\N}{\mathbb{N}}
\newcommand{\E}{\mathbb{E}}
\newcommand{\h}{\mathcal{H}}

\newcommand{\pr}{\mathbb{P}}
\newcommand{\1}{\mathbf{1}}

\newcommand{\da}{\downarrow}

\newcommand{\D}{{\rm d}}

\newcommand{\bd}{\begin{definition}}
\newcommand{\ed}{\end{definition}}
\newcommand{\bt}{\begin{theorem}}
\newcommand{\et}{\end{theorem}}
\newcommand{\bp}{\begin{proposition}}
\newcommand{\ep}{\end{proposition}}
\newcommand{\br}{\begin{remark}}
\newcommand{\er}{\end{remark}}
\newcommand{\bl}{\begin{lemma}}
\newcommand{\el}{\end{lemma}}
\newcommand{\bc}{\begin{corollary}}
\newcommand{\ec}{\end{corollary}}
\newcommand{\beq}{\begin{equation}}
\newcommand{\eeq}{\end{equation}}

\begin{document}
\baselineskip=13pt

\title{Mean density of inhomogeneous Boolean models with\\ lower dimensional
 typical grain}
\author{Elena Villa}
\date{\small Dept. of Mathematics, University of Milan, via Saldini 50,
20133 Milano, Italy\\
email: elena.villa@mat.unimi.it}

\maketitle

\begin{abstract}
The mean density of a random closed set $\Theta$ in $\R^d$ with
Hausdorff dimension $n$ is the Radon-Nikodym derivative of the
expected measure $\E[\h^n(\Theta\cap\cdot\,)]$ induced by $\Theta$
with respect to the usual $d$-dimensional Lebesgue measure. We
consider here inhomogeneous Boolean  models with lower dimensional
typical grain. Under general regularity assumptions on the typical
grain, related to the existence of its Minkowski content, and on
the intensity measure of the underlying Poisson point process, we
prove an explicit formula for the mean density. The proof of such
formula provides as by-product   estimators for the mean density
 in terms of the empirical capacity functional,
which turns to be closely related to the well known random
variable density estimation by histograms in the extreme case
$n=0$. Particular cases and examples are also discussed.
\end{abstract}


 AMS Classification 2000: 60D05; 28A75; 60G55; 49Q15

\section{Introduction}
In many real applications it is of interest to study random closed
sets at different Hausdorff dimensions and their induced random
measure; in particular
   several problems are related to the
estimation of the density, said \emph{mean density}, of the
\emph{expected measure}
$$\E[\mu_\Theta](B):=\mathbb{E}[\mathcal{H}^n(\Theta\cap B)]
\qquad \forall B\in\mathcal{B}_{\R^d}$$ induced by a random closed
set
$\Theta:(\Omega,\mathfrak{F},\pr)\to(\mathbb{F},\sigma_\mathbb{F})$
in $\R^d$ with Hausdorff dimension $n$ ($\h^n,\, {\cal
B}_{\R^d},\, \mathbb{F}$ and  $\sigma_\mathbb{F}$ denote here  the
$n$-dimensional Hausdorff measure, the Borel $\sigma$-algebra of
$\R^d$, the class of the closed subsets in $\mathbb{R}^d$ and the
$\sigma$-algebra generated by the so-called hit-or-miss topology
\cite{matheron}, respectively).
 While the extreme cases $n=d$ and
$n=0$ are easy to handle with elementary analytical tools,
problems arise when $0<n<d$. Indeed, while the expected measure
induced by a $d$-dimensional random set is always absolutely
continuous with respect to the $d$-dimensional Lebesgue measure
and its density can be easily obtained by applying Fubini's
theorem (as stated in the early  work by Robbins \cite{robbins}),
and the mean density of a random point is given by its probability
density function (and so the problem of its estimation has been
largely solved since long in nowadays standard literature, by
means of either histograms, or kernel estimators (e.g. see
\cite{silverman})),  when we deal with a general lower dimensional
random closed set it can be more demanding to check that the
induced expected measure is absolutely continuous and to compute
and estimate its density. Even if inhomogeneous random closed sets
appear frequently in real applications, only the stationary case
has been extensively studied so far (e.g., see
\cite{Benes-Rataj,SKM} and reference therein).
 Indeed,
although some results about the mean densities of certain
inhomogeneous random closed sets are available in current
literature  (e.g., see \cite{benes et
al,lars,HLW-contact,schneider,weil}), mainly via tools from
integral geometry and stereology, the pointwise estimation of the
mean density of nontrivial lower dimensional inhomogeneous random
closed sets is still an open problem, whenever gradient structures
\cite{HMS} or local stationarity cannot be assumed.

 Our goal is to obtain a
pointwise result for the mean density
 of a wide class of
inhomogeneous random sets which might be of interest from a
statistical point of view as well. We consider here inhomogeneous
Boolean models with lower dimensional typical grain, being
considered basic random set models in stochastic geometry
\cite{CIME}, but we do not exclude that our approach could be
applied to other kinds of inhomogeneous lower dimensional random
closed sets  in further developments. The main result of the paper
is stated in Theorem~\ref{teo dens}, where we prove an explicit
formula for the mean density; we point out that, even if such
formula could be obtained in a more direct way via the well known
Campbell's formula (see Remark~\ref{remark Campbell}), the proof
we propose here provides as a by-product estimators of the mean
density, answering, in the case of Boolean models,  to an open
problem in \cite{AKV}. To this end,
 the classical
Minkowski content notion plays a central role throughout the
paper; as a matter of fact, the regularity assumptions on the
typical grain we require are closely related to the well known
general assumptions which guarantee the existence of the Minkowski
content of deterministic closed sets.

 The paper is organized as
follows. In Section~\ref{prelim} we recall some basic definitions
and preliminary results useful for the sequel. In Section~\ref{sec
Minkowski} we prove a generalization of the Minkowski content of
closed subsets of $\R^d$ (Theorem~\ref{Mink gen}); we will apply
such result in the final part of the proof of Theorem~\ref{teo
dens}. In Section~\ref{sec mean density} we introduce a wide class
of Boolean models in $\R^d$ with lower dimensional typical grain.
We prove that the expected measures induced by the Boolean models
 we consider are absolutely continuous with respect to the usual
Lebesgue measure (Lemma~\ref{lemma a.c.}) and we provide an
explicit formula for their density. Such formula simplifies in the
special cases in which the Boolean model is assumed to be
stationary (Corollary~\ref{cor stat}), or to have deterministic
typical grain (Corollary~\ref{cor deterministic}). A simple
example of inhomogeneous segment Boolean model and links with
current literature are also provided. Finally, the problem of the
estimation of the mean density is considered in Section~\ref{sec
estimation}; in particular, we define an estimator of the mean
density (Proposition~\ref{prop stimatore}), which can be
considered as the generalization to the case $0<n<d$ of the
classical random variable density estimation by histograms in the
extreme case $n=0$.

\section{Basic notation and preliminaries}\label{prelim}
 We recall
that $\mathcal{H}^0$ is the usual counting measure,
$\mathcal{H}^n(B)$ coincides with the classical $n$-dimensional
measure of $B$  for $1\leq n<d$ integer if
$B\in\mathcal{B}_{\mathbb{R}^d}$  is contained in a $C^1$
$n$-dimensional manifold embedded in $\mathbb{R}^d$,
$\mathcal{H}^d(B)$ coincides with the usual $d$-dimensional
Lebesgue measure of $B$ for any Borel set $B\subset\mathbb{R}^d$.
Throughout the paper $\D x$ stands for $\h^d(\D x)$.  A closed
subset $S$ of $\mathbb{R}^d$ is said to be \emph{countably
$\mathcal{H}^{n}$-rectifiable} if there exist countably many
$n$-dimensional Lipschitz graphs $\Gamma_i\subset\mathbb{R}^d$
such that $S\setminus\cup_i\Gamma_i$ is
$\mathcal{H}^n$-negligible. (For definitions and basic properties
of Hausdorff measure and rectifiable sets see, e.g.,
\cite{AFP,Falconer,Fed}.) We  call \emph{Radon measure} in
$\mathbb{R}^d$ any nonnegative and $\sigma$-additive set function
$\mu$ defined on $\mathcal{B}_{\mathbb{R}^d}$ which is finite on
bounded sets, and we write $\mu\ll\h^n$ to say that $\mu$ is
absolutely continuous with respect to $\h^n$.
 We
will say that a random closed set $\Theta$ satisfies a certain
property (e.g. $\Theta$ has Hausdorff dimension $n$) if
$\Theta(\omega)$ satisfies that property for
  $\pr$-a.e.~$\omega\in\Omega$.
\\
  Throughout the paper we will
  consider countably $\h^n$-rectifiable random closed sets in
  $\R^d$, with $1\leq n\leq d-1$ integer, such that $\E[\mu_\Theta]$ is a Radon measure;
   the particular cases
  $n=0$ and $n=d$ are trivial.  (For a discussion about measurability
of $ \mathcal{H}^n(\Theta)$ we refer to \cite{Zahle,MB}.)
\\Whenever $\E[\mu_\Theta]$ is  absolutely continuous with respect to $\h^d$, the following
definition is given \cite{IAS,PhD}
\begin{definition}[Absolute continuity in mean and mean density]
Let $\Theta$ be a countably $\h^n$-rectifiable random closed set
in $\mathbb{R}^d$ such that $\E[\mu_\Theta]$ is a Radon measure.
We say that $\Theta$ is \emph{absolutely continuous in mean} if
$\mathbb{E}[\mu_{\Theta}]\ll\h^d$. In this case we call \emph{mean
density} of $\Theta$, and denote by $\lambda_{\Theta}$, the
Radon-Nikodym derivative of $\mathbb{E}[\mu_{\Theta}]$ with
respect to $\h^d$.
\end{definition}
\br\label{rem n=0 n=d} {\rm In the case $n=0$ with $\Theta=X$
random point in $\R^d$, we have that $\Theta$ is absolutely
continuous in mean if and only if $X$ admits a probability density
function $f_X$, and so $\lambda_\Theta=f_X$.\\
 On the other hand, it is easy to see
(as an application of Fubini's theorem in $\Omega\times\R^d$ with
the product measure $\pr\times\h^d$) that any $d$-dimensional
random closed set $\Theta$ in $\R^d$ is absolutely continuous in
mean with mean density $\lambda_\Theta(x)=\pr(x\in\Theta)$ for
$\h^d$-a.e~$x\in\R^d$. }\er

 The problem of the approximation of
the mean densities in the general setting of spatially
inhomogeneous processes has been recently faced in \cite{AKV},
where an approximation, in weak form, of the mean density for
sufficiently regular random closed sets is given in terms of their
$d$-dimensional enlargement by Minkowski addition. More precisely,
denoting by $ S_{\oplus r}:=\left\{x\in\mathbb{R}^d:\
\text{$\exists y\in S$ with $|x-y|\leq r$}\right\}$ the closed
$r$-neighborhood of a closed set $S\subset\R^d$, and by $B_r(x)$
the closed ball centered in $x$ with radius $r$, for any compact
window $W\subset\R^d$ let
$\Gamma_W(\Theta):\Omega\longrightarrow\R$ be the function so
defined:
\begin{multline}\Gamma_W(\Theta):=\sup\bigl\{\gamma\geq
0\,:\, \exists \mbox{ a probability measure }\eta\ll\mathcal{H}^n
\mbox{ such that }
\label{cond gamma}\\
\eta(B_r(x))\geq\gamma r^n\quad\forall x\in\Theta\cap W_{\oplus
1},\ r\in(0,1)\bigr\};
\end{multline}
then the following theorem holds  \cite{AKV}.

\begin{theorem}\label{teorema riassuntivo}
Let $\Theta$ be a countably $\h^n$-rectifiable random closed set
in $\mathbb{R}^d$ such that $\E[\mu_\Theta]$ is a Radon measure.
 Assume
that for any compact window $W\subset\mathbb{R}^d$ there exists a
random variable $Y$ with $\mathbb{E}[Y]<\infty$, such that
$1/\Gamma_W(\Theta)\leq Y$ almost surely.
 If $\Theta$ is absolutely continuous in mean, then
\begin{equation}\label{main2}
\lim_{r\da 0}\int_A \frac{\mathbb{P}(x\in {\Theta}_{\oplus
r})}{b_{d-n}r^{d-n}} {\rm d} x= \int_A \lambda_{\Theta}(x) {\rm d}
x
\end{equation}
for any bounded Borel set $A\subset\mathbb{R}^d$ with
$\h^d(\partial A)=0$.
\end{theorem}
 Hence the above theorem gives a \emph{weak} result for the mean
density of very general lower dimensional random closed sets. In
order to obtain a \emph{pointwise} result, using  the fact that
$A$ is arbitrary, we should prove that limit and integral in
\eqref{main2} can be exchanged; in such a way we could state that
\beq\lambda_{\Theta}(x)=\lim_{r\da 0}\frac{\mathbb{P}(x\in
{\Theta}_{\oplus r})}{b_{d-n}r^{d-n}} \qquad \mbox{for
}\h^d\mbox{-a.e.}~x\in\R^d,\label{scambio}\eeq and, as by-product,
the right side of the above equation could suggest estimators for
$\lambda_\Theta(x)$ in terms of the capacity functional of
$\Theta$.
 The proof of the
validity of this formula for absolutely continuous in mean random
sets seems to be a quite delicate problem, with the only exception
of the stationary ones and  the extreme cases $n=d$ and $n=0$. In
Theorem~\ref{teo dens} we prove that \eqref{scambio} holds for a
wide class of (inhomogeneous) \emph{Boolean models} and, in
particular, we give an explicit representation of the mean density
$\lambda_\Theta$ in terms of the \emph{intensity} and of the
\emph{typical grain} of the process. To this end we  recall now
the  definition  (typically standard) of Boolean model and we
refer to \cite{CIME,Daley} for details about point processes and
related concepts and results. \bd[Boolean model] Let
$\Psi=\{x_i\}_{i\in\N}$ be Poisson point process in $\R^d$ with
intensity $f$ and let $\{Z_i\}_{i\in\N}$ be a sequence of i.i.d.
random compact sets in $\R^d$, which are also independent of the
Poisson process $\Psi$. Denoting by $Z_0$ a further random compact
set of the same distribution as the $Z_i$'s and independent of
both of them and of $\Phi$, the resulting random set
$$\Theta:=\bigcup_i (x_i+Z_i)
$$ is said \emph{(inhomogeneous) Boolean model with intensity $f$
and typical grain} $Z_0$. \ed In our assumptions $Z_0$ will be a
lower dimensional random closed set uniquely determined by a
random quantity in a suitable \emph{mark space} $\bf K$, so that
$Z_0(s)$ is a compact subset of $\R^d$ containing the origin for
any $s\in\bf K$. We remind that in common literature is usually
assumed that \beq\E[{\rm card}\{i\,:\, (x_i+Z_i)\cap
K\neq\emptyset\}]<\infty \qquad \forall \mbox{ compact
}K\subset\R^d\label{cond Bool1}\eeq (card stands for cardinality),
and that
 a Boolean model as above  can
 be described by a  Poisson point process in $\R^d\times\bf K$
 with intensity measure $\Lambda(\D(x,s))=f(x)\D xQ(\D s)$, where
 $Q$ is a probability measure on $\bf K$ representing the distribution
 of the typical grain.

\section{A generalization of the Minkowski content of  closed
sets}\label{sec Minkowski} Denoted by $b_k$ the volume
 of the unit ball in $\R^k$,
the $n$-dimensional \emph{Minkowski content} ${\cal M}^n(S)$ of a
closed set $S\subset\mathbb{R}^d$ is defined by
$$
{\cal M}^n(S):=\lim_{r\da 0}\frac{\h^d(S_{\oplus
r})}{b_{d-n}r^{d-n}}
$$
whenever the limit exists finite.
\\
General results about the existence of the Minkowski content of
closed subsets in $\mathbb{R}^d$ are known in literature,
 related to rectifiability properties of the involved
sets. In particular, the following theorem is proved in \cite{AFP}
(p.\,110). \begin{theorem}\label{teo mink cont}  Let
$S\subset\mathbb{R}^d$ be a countably
$\mathcal{H}^{n}$-rectifiable compact set and assume that
\begin{equation}\eta(B_r(x))\geq \gamma r^n \qquad \forall x\in S,\,\,\forall r\in(0,1)
\label{condizione eta}\end{equation} holds for some $\gamma>0$ and
some Radon measure $\eta$ in $\mathbb{R}^d$ absolutely continuous
with respect to $\mathcal{H}^n$. Then
$\mathcal{M}^n(S)=\mathcal{H}^n(S)$. \end{theorem} Note that such
theorem extends the well-known Federer's result (\cite{Fed},
p.\,275) about the existence of ${\cal M}^n(S)$ for
$n$-rectifiable compact sets $S\subset\R^d$ (i.e. $S$ is
representable as the image of a compact subset of $\mathbb{R}^n$
by a Lipschitz function from $\mathbb{R}^n$ to $\mathbb{R}^d$) to
countably $\mathcal{H}^{n}$-rectifiable compact sets (see Remark
2.3 in \cite{acolev}). Moreover, in many applications condition
\eqref{condizione eta} is satisfied with
$\eta(\cdot)=\h^n(\widetilde{S}\cap\cdot\,)$ for some closed set
$\widetilde{S}\supseteq S$, and  for $n=d-1$, it is not hard to
check that such condition is satisfied by all sets with Lipschitz
boundary (see \cite{acolev,AFP}).
\\
We state now the main result of this section. \bt\label{Mink gen}
Let $\mu$ be a positive measure in $\R^d$ absolutely continuous
with respect to $\h^d$ with density $f$ such that
\begin{itemize}
\item[i)] $f$ is locally bounded (i.e.~$\sup_{x\in K}f(x)<\infty$
for any compact $K\subset\R^d$); \item[ii)] the set of all
discontinuity points of $f$ is $\h^n$-negligible.
\end{itemize}
Let $S\subset\mathbb{R}^d$ be a countably
$\mathcal{H}^{n}$-rectifiable compact set such that condition
\eqref{condizione eta} holds for some $\gamma >0$ and some
probability measure $\eta$ in $\R^d$ absolutely continuous with
respect to $\h^n$. Then

 $$\lim_{r\da 0}\frac{\mu(S_{\oplus r})}{b_{d-n}r^{d-n}}=\int_S f(x)\h^n(\D
 x).$$
 \et

\br\label{remark finite}{\rm The above theorem may be seen as a
generalization of Theorem~\ref{teo mink cont}; indeed, the
classical Minkowski content follows as particular case by choosing
$f\equiv 1$ and noticing that if  a Radon measure $\eta$ as in
Theorem~\ref{teo mink cont} exists, then
 it can be assumed to be a probability
measure without loss of generality. Indeed, it is sufficient to
consider the measure
$\tilde\eta(\cdot):=\eta(W\cap\cdot\,)/\eta(W)$, where $W$ is a
compact subset of $\R^d$ such that $S_{\oplus 1}\subset W$; it is
clear that $\tilde\eta$ is a probability measure satisfying
$$\tilde{\eta}(B_r(x))\geq \frac{\gamma}{\eta(W)} r^n \qquad \forall x\in S,\,\,\forall r\in(0,1).
 $$

 Furthermore, a classical result from geometric measure theory (e.g., see \cite{AFP} Theorem 2.56)
 tells us that if $\mu$ is a positive Radon measure on $\R^d$ and $B\in\mathcal{B}_{\R^d}$ such that
 $$
\limsup_{r\downarrow 0}\frac{\mu(B_{r}(x))}{b_nr^n}\geq
t\in(0,\infty)\qquad\forall x\in B, $$ then  $\mu(\cdot)\geq
t\h^n(B\cap \cdot)$. Hence, any set $S\subset\R^d$ as in Theorem
\ref{Mink gen} has  finite $\h^n$-measure.}\er In order to make
the proof of Theorem~\ref{Mink gen} more readable, we remind that
 Lemma 15 and
Lemma 16 in \cite{AKV} provide a local version of Theorem~\ref{teo
mink cont} and an upper bound for the Minkowski content of compact
sets in $\R^d$, respectively; for our purpose we summarize as
follows. \bl\label{lemma intersezione} If $S\subset\mathbb{R}^d$
is a countably $\mathcal{H}^{n}$-rectifiable compact set such that
condition \eqref{condizione eta} holds for some $\gamma >0$ and
some finite measure $\eta$ in $\R^d$ absolutely continuous with
respect to $\h^n$,
 then \beq\frac{\h^d(S_{\oplus
r})}{b_{d-n}r^{d-n}}\leq \frac{\eta(\R^d)}{\gamma}2^n
4^d\frac{b_d}{b_{d-n}}\qquad \forall r<2,\label{lemma AKV}\eeq and
 \begin{equation}\lim_{r\da 0}\frac{\h^d(S_{\oplus r}\cap
A)}{b_{d-n}r^{d-n}}=\mathcal{H}^n(S\cap A)\label{tesi} \eeq for
any $A\in\mathcal{B}_{\mathbb{R}^d}$ such that $
\mathcal{H}^n(S\cap\partial A)=0$. \el

\vspace{3mm}\noindent \emph{Proof of Theorem~\ref{Mink gen}}. \ It
is well known that, since $f\geq 0$, there exists an increasing
sequence $\{f_k\}_{k\in\N}$ of step functions
$$f_k(x)=\sum_{j=1}^{N(k)}a_j^{(k)}\1_{A_j^{(k)}}(x),\quad a_j^{(k)}\geq 0,\  N(k)\in\N,$$ converging to
$f$. \\
By ii), such sequence $\{f_k\}$ can be chosen so that
$\h^n(S\cap\partial
A_j^{(k)})=0$ for all $j,k$.\\
Let
$$g_k(r):=\frac{\int_{S_{\oplus r}}f_k(x)\D x}{b_{d-n}r^{d-n}}. $$
We may observe that\\
(a)\quad $\displaystyle\lim_{r\da
0}g_k(r)=\sum_{j=1}^{N(k)}a_j^{(k)}\lim_{r\da
0}\frac{\h^d(S_{\oplus r}\cap
A_j^{(k)})}{b_{d-n}r^{d-n}}\stackrel{\eqref{tesi}}{=}\sum_{j=1}^{N(k)}a_j^{(k)}\h^n(S\cap
A_j^{(k)})= \int_Sf_k(x)\h^n(\D x).$ \\Since $S$ is compact with
$\h^n(S)<\infty$ by Remark~\ref{remark finite}, and $f$ is locally
bounded, it follows that
$\displaystyle\int_Sf_k(x)\h^n(\D x)<\infty$.\\[2mm]
(b)\quad $\displaystyle\lim_{k
\to\infty}g_k(r)=\frac{\int_{S_{\oplus r}}f(x)\D
x}{b_{d-n}r^{d-n}}\quad \mbox{uniformly in }(0,1)$; indeed:\\
we have that  $f_k\uparrow f$ uniformly in $S_{\oplus r}$ for any
$r>0$ because $f$ is bounded in $S_{\oplus r}$, being $S_{\oplus
r}$  compact,
 i.e.~for all
$\varepsilon>0$ there exists $k_0$ such that \beq \sup_{x\in
S_{\oplus r}}|f_k(x)-f(x)|<\varepsilon\qquad \forall
k>k_0.\label{unif f}\eeq Hence, for any $\varepsilon>0$, for all
$k>k_0$
$$\left|g_k(r)-\frac{\int_{S_{\oplus r}}f(x)\D
x}{b_{d-n}r^{d-n}}\right|\leq\frac{\int_{S_{\oplus
r}}|f_k(x)-f(x)|\D x}{b_{d-n}r^{d-n}} \stackrel{\eqref{unif f}}{<}
\varepsilon\frac{\h^d(S_{\oplus r})}{b_{d-n}r^{d-n}}
\stackrel{\eqref{lemma AKV}}{\leq}\varepsilon\frac{1}{\gamma}2^n
4^d\frac{b_d}{b_{d-n}}\qquad \forall r\in(0,1).
$$

\vspace{3mm} \noindent As a consequence of (a) and (b) $\lim_{r\da
0}\lim_{k\to\infty}g_k(r)=\lim_{k\to\infty}\lim_{r\da 0}g_k(r)$.
Thus the following chain of equalities holds
\begin{multline*} \lim_{r\da 0}\frac{\mu(S_{\oplus
r})}{b_{d-n}r^{d-n}}=\lim_{r\da 0}\frac{\int_{S_{\oplus r}}f(x)\D
x}{b_{d-n}r^{d-n}}=\lim_{r\da 0}\frac{\int_{S_{\oplus
r}}\lim_{k\to\infty}f_k(x)\,\D x}{b_{d-n}r^{d-n}} =\lim_{r\da
0}\lim_{k\to\infty}\frac{\int_{S_{\oplus r}}f_k(x)\D x}{b_{d-n}r^{d-n}}\\
=\lim_{k\to\infty}\lim_{r\da 0}\frac{\int_{S_{\oplus r}}f_k(x)\D
x}{b_{d-n}r^{d-n}} \stackrel{({\rm a})}{=}\lim_{k\to\infty} \int_S
f_k(x)\h^n(\D x)=\int_Sf(x)\h^n(\D x).
\end{multline*}
{\hfill $\square$\vspace{2mm}}

 \section{Mean density of inhomogeneous Boolean models}\label{sec mean density}
 Let $\Theta$ the Boolean model so
defined
$$\Theta(\omega):=\bigcup_{(x_i,s_i)\in\Phi(\omega)}x_i+Z_0(s_i) \qquad \forall \omega\in\Omega, $$
with $\Phi$ Poisson process in $\R^d\times{\bf K}$ and
 $Z_0$ typical grain in $\R^d$, satisfying the usual condition
\eqref{cond Bool1} for Boolean models.  To lighten the notations,
from now on we denote by\\
--\,  $Z^{x,s}:=x-Z_0(s)$ $\forall (x,s)\in\R^d\times\bf K$;
\\
--\, $\delta$ the random variable on $\bf K$
 so defined $\delta(s):={\rm diam}Z_0(s)$, where diam stands for
 diameter;
\\
--\,  $\E_Q$ the expectation with respect to the probability
measure $Q$ on $\bf K$.

Note that the condition \eqref{cond Bool1} is equivalent to
 $$\E[\Phi(\{(y,s)\in\R^d\times{\bf K}\,:\,
(y+Z_0(s))\cap B_R(0)\neq\emptyset\})]<\infty \qquad \forall
R>0,$$ and so, in terms of the intensity measure $\Lambda$ of
$\Phi$, to
 \beq\int_{\bf
K}\int_{(-Z_0(s))_{\oplus R}}\Lambda(\D y\times \D s) <\infty
\qquad\forall R>0, \label{cond Bool 2}\eeq where $-Z_0$ is the
symmetric of $Z_0$ with respect to the origin.

 {\bf Assumptions:}
Let $\Phi$ have intensity measure $\Lambda(\D y\times\D s)=f(y)\D
y Q(\D s)$ satisfying \eqref{cond Bool 2}; further, let us assume
that the following conditions
on $f$ and $Z_0$ are fulfilled:\\[1mm]
{\bf (A1)} \ $Z_0(s)$ is a countably $\h^n$-rectifiable
 compact set in $\R^d$ for $Q$-a.e.~$s\in\bf K$. Further there exist $\gamma>0$ and a random closed set $\widetilde{Z}_0\supseteq Z_0$
 with
 $\E_Q[\h^n(\widetilde{Z}_0)]<\infty$
 such that, for $Q$-a.e.~$s\in\bf
 K$,  \beq\h^n(\widetilde{Z}_0(s)\cap B_r(x))\geq \gamma r^n
 \qquad \forall x\in Z_0(s),\,\,\forall r\in(0,1). \label{cond Z
 tilde}\eeq
{\bf(A2)} \   the set of all discontinuity points of $f$ is
$\h^n$-negligible and $f$ is locally bounded  such that for any
compact set $K\subset\R^d$ \beq\sup_{y\in K_{\oplus
\delta}}f(y)\leq \xi_K\label{sup}\eeq holds for some  random
variable $\xi_K$ with $\E_Q[\h^n(\widetilde{Z}_0)\xi_K]<\infty$.
\bt[Main result]\label{teo dens} Any Boolean model $\Theta$ as in
Assumptions is absolutely continuous in mean with mean density
$\lambda_\Theta$ given by \beq\lambda_\Theta(x)=\int_{\bf
K}\int_{Z^{x,s}}f(y)\h^n(\D y)\,Q(\D s)\qquad \mbox{for
}\h^d\mbox{-a.e.}~x\in\R^d. \label{lambda}\eeq \et

 Before proving the above theorem, let us make a few
observations about the above Assumptions in order to clarify their
level of generality. \begin{itemize} \item
 Condition \eqref{sup} is trivially satisfied whenever  $f$ is
bounded, or $f$ is locally bounded and ${\rm diam}Z_0\leq
c\in\R_+$ $Q$-almost surely. \item Assumption (A1) is satisfied by
a great deal of typical grains $Z_0$ with $\widetilde{Z}_0=Z_0$ or
$\widetilde{Z}_0=Z_0\cup \widetilde{A}$, with $\widetilde{A}$
sufficiently regular random closed set to control the case when
$Z_0$ can be arbitrarily small (see \cite{AKV,PhD}), and it
assures that $Z_0$ is a countably $\h^n$-rectifiable random
compact set with finite expected $\h^n$-measure.
 \item   {\bf Stationary case.} If $\Theta$ is stationary with
  $f\equiv c>0$, then only the regularity assumption (A1) on the typical grain $Z_0$
   is required. Indeed, the assumption (A2) is trivial and the usual condition \eqref{cond Bool
   2}, which in this case becomes
   $$\E_Q[\h^d((Z_0)_{\oplus
R})]<\infty\qquad \forall R>0, $$  is satisfied thanks to the
assumption (A1): for all $R<2$ by Lemma~\ref{lemma intersezione}
with $\eta(\cdot)=\h^n(\widetilde{Z}_0\cap\cdot\,)$ we have that
$\E_Q[\h^d((Z_0)_{\oplus R})]\leq \E[\h^n(\widetilde{Z}_0)]2^n4^d
b_dR^{d-n}/\gamma$, while for all $R\geq 2$, by repeating  the
same argument of the proof of the quoted proposition, it is easy
to check that $\E_Q[\h^d((Z_0)_{\oplus R})]\leq
\E[\h^n(\widetilde{Z}_0)]2^n(4R)^d b_d/\gamma$. \item {\bf
Deterministic typical grain.} If $\Theta$ has deterministic
typical grain $Z_0\subset\R^d$ the above Assumptions can be
replaced by: $\Theta$ has intensity $f$ and typical grain $Z_0=S$
as in the hypotheses of Theorem~\ref{Mink gen}. (See
Corollary~\ref{cor deterministic} for details.)
\end{itemize}

\begin{lemma}\label{lemma a.c.}
For any Boolean model $\Theta$ as in Assumptions, $\E[\mu_\Theta]$
is a locally finite measure absolutely continuous with respect to
$\h^d$.
\end{lemma}
\begin{proof} For any $R>0$ let ${\cal B}^R:=\{(y,s)\in\R^d\times{\bf K}\,:\,
(y+Z_0(s))\cap B_R(0)\neq\emptyset\}$; then
$$
\E[\h^n(\Theta\cap B_R(0))]=\E[\E[\h^n(\Theta\cap
B_R(0))\,|\,\Phi({\cal B}^R)]]\leq \E_Q[\h^n(Z_0)]\E[\Phi({\cal
B}^R)]<\infty
$$
by \eqref{cond Bool 2} and condition (A1); so $\E[\mu_{\Theta}]$
is locally finite.

 By contradiction, let
  $\mathbb{E}[\mu_\Theta]$ be not absolutely continuous with respect to $\h^d$;
then there exists
 $A\subset\mathbb{R}^d$ with
$\h^d(A)=0$ such that $\mathbb{E}[\h^n(\Theta\cap A)]>0$. In
particular, $$ 0<\pr(\h^n(\Theta\cap
A)>0)\leq\pr\Big(\sum_{(x_i,s_i)\in\Phi}\h^n((x_i+Z_0(s_i))\cap
A)>0\Big)=\pr(\Phi({\cal A})>0),
$$
 where
 $$\mathcal{A}:=\{(y,s)\in\R^d\times{\bf K} \,:\, \mathcal{H}^n((y+Z_0(s))\cap A)>0)\}. $$
Denoting by $\mathcal{A}_s:=\{y\in\R^d \,:\, (y,s)\in\mathcal A\}$
the section of $\mathcal{A}$ at $s\in\bf K$, we may apply Fubini's
theorem to get
$$\int_{{\cal A}_s}\h^n((y+Z_0(s))\cap A)\D y=\int_{{\cal A}_s}\Big(\int_{Z_0(s)}\1_{A-y}(x)\h^n(\D x)
\Big)\D y = \int_{Z_0(s)}\Big( \int_{{\cal A}_s}\1_{A-x}(y)\D
y\Big)\h^n(\D x)=0,
$$
because $\h^d(A)=0$. Being the function
$y\mapsto\mathcal{H}^n((y+Z_0(s))\cap A)$ strictly positive in
${\cal A}_s$, we conclude that  $\h^d(\mathcal{A}_s)=0$ for all
$s\in{\bf K}$. Then it follows
$$\mathbb{E}[\Phi(\mathcal{A})]=\int_\mathcal{A}\Lambda(\D y\times \D s)=\int_{\bf K}\Big(\int_{{\cal A}_s}
f(y)\D y\Big) Q(\D s)=0;$$ but this is impossible, because
$\mathbb{P}(\Phi(\mathcal{A})>0)>0$ implies
$\mathbb{E}[\Phi(\mathcal{A})]>0. $
\end{proof}

\vspace{3mm}\noindent \emph{Proof of Theorem~\ref{teo dens}}. \
Clearly $\Theta$ is a countably $\h^n$-rectifiable random closed
set in $\R^d$, and by Lemma~\ref{lemma a.c.} it is absolutely
continuous in mean, so $\E[\mu_\Theta]=\lambda_\Theta\h^d$ for
some integrable function
$\lambda_\Theta$ on $\R^d$.\\
Let $W$ be a fixed compact subset of $\R^d$. For any
$\omega\in\Omega$ let us consider the set
$$\widetilde{\Theta}(\omega):=\bigcup_{(x_i,s_i)\in\Phi(\omega)}x_i+\widetilde{Z}_0(s_i)$$
and the probability measure $\eta_W$ absolutely continuous with
respect to $\h^n$ so defined
$$\eta_W(B):=\frac{\h^n(\widetilde{\Theta}(\omega)\cap W_{\oplus 1}\cap B)}
{\h^n(\widetilde{\Theta}(\omega)\cap W_{\oplus 1})}\qquad \forall
B\in\mathcal{B}_{\R^d}. $$
 Note that
for any $x\in\Theta(\omega)\cap W_{\oplus 1}$ there exists
$(\bar{x},\bar{s})\in\Phi(\omega)$ such that
$x\in\bar{x}+{Z}_0(\bar{s})$,
 and so
$$\eta_W(B_r(x))\geq\frac{\h^n(\widetilde{Z}_0(\bar{s})\cap B_r(x-\bar{x}))}
{\h^n(\widetilde{\Theta}(\omega)\cap W_{\oplus
1})}\stackrel{\eqref{cond Z
tilde}}{\geq}\frac{\gamma}{\h^n(\widetilde{\Theta}(\omega)\cap
W_{\oplus 1})}r^n\qquad \forall r\in(0,1).
$$ Thus, here  the function $\Gamma_W(\Theta)$ defined
in \eqref{cond gamma} is such that
$\Gamma_W(\Theta)\geq\gamma/\h^n(\widetilde{\Theta}\cap W_{\oplus
1})$; by the same argument in the proof of Lemma~\ref{lemma a.c.}
it is easy to check that $\E[\h^n(\widetilde{\Theta}\cap W_{\oplus
1})]<\infty$, then Theorem~\ref{teorema riassuntivo} applies and
we get \beq\E[\mu_\Theta](A)=\lim_{r\da 0}\int_A
\frac{\mathbb{P}(x\in {\Theta}_{\oplus r})}{b_{d-n}r^{d-n}} {\rm
d} x \label{E=lim}\eeq for any bounded set
$A\in\mathcal{B}_{\mathbb{R}^d}$ such that $
\h^d(\partial A)=0. $ \\
Let us denote by ${\cal Z}^{x,r}$ the subset of $\R^d\times{\bf
K}$ so defined
$${\cal Z}^{x,r}:= \{(y,s)\in\R^d\times{\bf K}\,:\,
x\in(y+Z_0(s))_{\oplus r} \}=\{(y,s)\in\R^d\times{\bf K}\,:\, y\in
Z^{x,s}_{\oplus r} \},
$$ and observe that $\forall x\in A,\,\,\forall r<2$
\begin{eqnarray}
\frac{\pr(x\in\Theta_{\oplus
r})}{b_{d-n}r^{d-n}}&=&\frac{\pr(\Phi({\cal
Z}^{x,r})>0)}{b_{d-n}r^{d-n}}=\frac{1-e^{-\Lambda({\cal
Z}^{x,r})}}{b_{d-n}r^{d-n}}\label{prob frac}\\
&\leq&\frac{\Lambda({\cal Z}^{x,r})}{b_{d-n}r^{d-n}}
=\frac{1}{b_{d-n}r^{d-n}}\int_{\bf
K}\int_{Z^{x,s}_{\oplus r}}f(y)\D y\,Q(\D s)\nonumber\\
&\leq&\int_{\bf K}\frac{ \h^d(Z^{x,s}_{\oplus
r})}{b_{d-n}r^{d-n}}\sup_{y\in A_{\oplus \delta(s)+2}}f(y)\, Q(\D
s)\nonumber\\
&\stackrel{\eqref{lemma AKV},\eqref{sup}}{\leq}&\frac{2^n 4^d
b_d}{\gamma b_{d-n}} \int_{\bf K}
\h^n(\widetilde{Z}_0(s))\xi_{A_{\oplus 2}}(s)Q(\D
s)\stackrel{(A2)}{<}\infty.\nonumber
\end{eqnarray}
 Then by the dominated convergence theorem we can exchange limit
and integral in Eq.~\eqref{E=lim}.

Similarly, we have that for all $r<2$
$$\frac{\int_{Z^{x,s}_{\oplus r}}f(y)\D y}{b_{d-n}r^{d-n}}\leq
\frac{\h^d((Z_0(s))_{\oplus r})}{b_{d-n}r^{d-n}}\sup_{y\in
Z^{x,s}_{\oplus r}} f(y)\leq\frac{2^n 4^d b_d}{\gamma
b_{d-n}}\h^n(\widetilde{Z}_0(s))\xi(s)$$ for some random variable
$\xi$ with $\E_Q[\h^n(\widetilde{Z}_0)\xi]<\infty$,  by
\eqref{lemma AKV} and (A2). Then, the dominated convergence
theorem implies that
\begin{equation}\lim_{r\da 0}\frac{\Lambda({\cal
Z}^{x,r})}{b_{d-n}r^{d-n}}=\int_{\bf K}\lim_{r\da
0}\frac{\int_{Z^{x,s}_{\oplus r}}f(y)\D y}{b_{d-n}r^{d-n}}\,Q(\D
s)=\int_{\bf K}\int_{Z^{x,s}}f(y)\h^n(\D y)\,Q(\D s),\label{tesi
th}\end{equation} where the last equality follows by applying
Theorem~\ref{Mink gen} with $\mu=f\h^d$ and $S=Z^{x,s}$.
\\ Summarizing, we have that for any bounded  set $A\in\mathcal{B}_{\mathbb{R}^d}$
with $\h^d(\partial A)=0 $
\begin{eqnarray*}\E[\mu_\Theta(A)]&=&\lim_{r\da 0}\int_A
\frac{\pr(x\in\Theta)}{b_{d-n}r^{d-n}}\D x=\int_A \lim_{r\da
0}\frac{\pr(x\in\Theta)}{b_{d-n}r^{d-n}}\,\D x\\
&\stackrel{\eqref{prob frac}}{=}&\int_A\lim_{r\da
0}\frac{1-e^{-\Lambda({\cal Z}^{x,r})}}{b_{d-n}r^{d-n}}\,\D x
=\int_A\lim_{r\da 0}\frac{\Lambda({\cal
Z}^{x,r})+o(r^{d-n})}{b_{d-n}r^{d-n}}\,\D
x\\&\stackrel{\eqref{tesi th}}{=}&\int_A\Big( \int_{\bf
K}\int_{Z^{x,s}}f(y)\h^n(\D y)\,Q(\D s) \Big)\D x.
\end{eqnarray*}
We conclude that $\E[\mu_\Theta]$ has density
$$\lambda_{\Theta}(x)=\int_{\bf K}\int_{Z^{x,s}}f(y)\h^n(\D y)\,Q(\D s)
 \qquad\mbox{for }\h^d\mbox{-a.e.~}x\in\R^d.$$
{\hfill $\square$\vspace{2mm}}

\br\label{remark Campbell} {\rm Formula \eqref{lambda} could also
be obtained in a more direct way by the well-known Campbell's
formula (e.g., see \cite{CIME}, p.\,28, and \cite{lars} for a
similar application), after having shown that
$\E[\mu_\Theta](\cdot)=\E[\sum_{i}\h^n((x_i+Z_i)\cap\,\cdot\,)]$
for any Boolean model $\Theta$ as in Assumptions, being zero the
probability that different grains overlap in a subset of $\R^d$ of
positive $\h^n$-measure. (To show such property, one could proceed
 arguing similarly to the final part of the proof of Lemma~\ref{lemma
a.c.}.) On the other hand, a  proof via Campbell's formula seems
not to give any hint to provide computable estimators for the mean
density (see Section~\ref{sec estimation} for a more detailed
discussion).

}\er

 In the particular cases in which $\Theta$ is stationary, or
the typical grain $Z_0$ is a deterministic subset of $\R^d$
satisfying the hypotheses of Theorem~\ref{teo mink cont} (with
$\eta$ probability measure by Remark~\ref{remark finite}),
Theorem~\ref{teo dens} specializes as follows.

\begin{corollary}[Stationary case]\label{cor stat}
Let $\Theta$ be a stationary Boolean model with intensity $f\equiv
c>0$ and typical grain $Z_0$ satisfying the assumption (A1). Then
$\Theta$ is absolutely continuous in mean with mean density
$$\lambda_{\Theta}(x)=c\E_Q[\h^n(Z_0)] \qquad \forall x\in\R^d.$$
\end{corollary}
\begin{proof}
We have observed that $\Theta$ satisfies the Assumptions and
$\lambda_\Theta$ is constant on $\R^d$. Then the thesis follows by
Theorem~\ref{teo dens}, noticing that $\h^n(Z^{x,s})=\h^n(Z_0(s))$
for all $(x,s)\in\R^d\times\bf K$ and that $\lambda_{\Theta}$ is
constant being $\E[\mu_\Theta]$ translation invariant.
\end{proof}

\begin{corollary}[Deterministic typical grain]\label{cor deterministic}
Let $\Theta$ be a Boolean model in $\R^d$ with deterministic
typical grain $Z_0\subset\R^d$. If the underlying Poisson process
$\Phi$ in $\R^d$ has intensity  $f$ locally bounded and such that
the set of all its discontinuity points is $\h^n$-negligible, and
if $Z_0$ is a countably $\mathcal{H}^{n}$-rectifiable compact set
such that condition \eqref{condizione eta} holds for some $\gamma
>0$ and some probability measure $\eta$ in $\R^d$ absolutely
continuous with respect to $\h^n$, then $\Theta$ is absolutely
continuous in mean with mean density
$$\lambda_{\Theta}(x)=
\int_{Z_0} f(x-y)\h^n(\D y) \qquad\mbox{for
}\h^d\mbox{-a.e.~}x\in\R^d.$$
\end{corollary}
\begin{proof}
Since $Z_0$ is compact  ${\rm diam}Z_0=\delta<\infty$, and we know
by Remark~\ref{remark finite} that $\h^n(Z_0)$ is finite. Then the
assumption (A2) and the condition \eqref{cond Bool 2} are easily
checked. It follows that
 $\E[\mu_\Theta]$
is locally finite and, by proceeding as in Lemma~\ref{lemma a.c.},
again we have that $\E[\mu_\Theta]\ll\h^d$.

 For any $y\in\R^d$ let $\eta^y$ be
 the measure on $\R^d$ so defined
$$\eta^y(B):=\eta(B-y)\qquad \forall B\in\mathcal{B}_{\R^d}.$$
Let $W$ be a fixed compact subset of $\R^d$ and  for any $\omega
\in\Omega$ consider the measure
$$\tilde\eta(B):=\frac{\sum_{x_i\in\Phi(\omega)}\eta^{x_i}(B)\1_{(x_i+Z_0)\cap W_{\oplus 1}\neq\emptyset}}
{{\rm card}\{x_i\in\Phi(\omega)\,:\, (x_i+Z_0)\cap W_{\oplus
1}\neq\emptyset\}}\qquad \forall B\in\mathcal{B}_{\R^d}.$$
 Note that
\\
 -- $\tilde\eta$ is a probability measure absolutely continuous with
respect to $\h^n$;
\\
-- for any $x\in\Theta(\omega)\cap W_{\oplus 1}$ there exists
$\bar{x}\in\Phi(\omega)$ such that $x\in\bar{x}+Z_0$,
 and so
$$\tilde\eta(B_r(x))\geq\frac{\gamma}{{\rm
card}\{x_i\in\Phi(\omega)\,:\, (x_i+Z_0)\cap W_{\oplus
1}\neq\emptyset\}}r^n\qquad \forall r\in(0,1).
$$ Then $\Theta$ satisfies the hypotheses of Theorem
\ref{teorema riassuntivo} with $Y={\rm
card}\{x_i\in\Phi(\omega)\,:\, (x_i+Z_0)\cap W_{\oplus
1}\neq\emptyset\}/\gamma$,  which has finite expectation thanks to
\eqref{cond Bool 2}, and so  again we have that Eq.~\eqref{E=lim}
holds for any bounded set $A\in\mathcal{B}_{\mathbb{R}^d}$ such
that $ \h^d(\partial A)=0 $. Finally, by proceeding as in the
final part of the proof of Theorem~\ref{teo dens} with
deterministic $Z_0$, we conclude that \eqref{scambio} still holds
and
$$\lambda_{\Theta}(x)=\int_{x-Z_0}f(y)\h^n(\D y)=
\int_{Z_0} f(x-y)\h^n(\D y) \qquad\mbox{for
}\h^d\mbox{-a.e.~}x\in\R^d.$$
\end{proof}

\br[The special case $n=d-1$]{\rm A well known tool in stochastic
geometry is the so-called {\em local spherical contact
distribution function}, defined as
$H_\Theta(r,x):=\pr(x\in\Theta_{\oplus r}\,|\,x\not\in\Theta)$ for
all $(r,x)\in\R_+\times \R^d$ (e.g, see \cite{HLW-contact}). If
$\Theta$ is a Boolean model as in Assumptions with
$(d-1)$-dimensional typical grain, then $\mathbb{P}(x\in\Theta)=0$
for $\mathcal{H}^d$-a.e.~$x\in \R^d$, and so the mean density
$\lambda_\Theta$ can be written in terms of $H_\Theta$ as well
$$\lambda_\Theta(x)\stackrel{\eqref{scambio}}{=}\frac{1}{2}\frac{\partial}{\partial
r}H_\Theta(r,x)_{|r=0}\qquad\mbox{for
}\h^d\mbox{-a.e.~}x\in\R^d.$$ Therefore known results about
$H_\Theta$ could be applied to $(d-1)$-dimensional Boolean models
to obtain further information on their mean density.  }\er

 We conclude
this section by giving a very simple example of inhomogeneous
Boolean model in $\R^d$, in order to show how the mean density can
be easily computed by applying Theorem~\ref{teo dens}.
\begin{example}[Segment Boolean model]{\rm  For sake of simplicity we consider a
Boolean model $\Theta$ of segments in
$\R^2$, but a similar example can be done in $\R^d$ with $d>2$.\\
So, let ${\bf K}=\R_+\times [0,2\pi]$
 and for all $s=(l,\alpha)\in\bf K$ let $Z_0(s)$ be the segment
 with length $l$ and orientation $\alpha$ so defined
$$Z_0(s):=\{(u,v)\in\R^2\,:\, u=\tau\cos\alpha,\ v=\tau\sin\alpha,
\ \tau\in[0,l] \}.$$  We consider the case in which both length
and orientation are random. Denoting by $L$ the $\R_+$-valued
random variable representing the length of $Z_0$ and by $\pr_L(\D
l)$ its probability law, let $\Phi$ be the marked Poisson process
in $\R^d\times\bf K$ having intensity measure $\Lambda(\D y\times
\D s)=f(y)\D yQ(\D s)$ with $f(u,v)=u^2+v^2$ and $Q(\D
s)=\frac{1}{2\pi}\D \alpha\pr_L(\D l)$ such that
$\int_{\R_+}l^3\pr_L(\D l)<\infty$. (This last assumption is to
guarantee that the usual condition \eqref{cond Bool1} is
satisfied; for a different intensity $f$ we might have a different
condition on the moments of $L$.) It is easily shown
\cite{AKV,PhD} that $\Theta$ satisfies the assumption (A1); by
noticing that
$$\int_{(-Z_0(s))_{\oplus R}}f(y)\D y\leq(l+R)^2(2lR+\pi R^2) \qquad\forall s=(l,\alpha)\in\bf K,$$
it follows that
$$
\int_{\bf K}\int_{(-Z_0(s))_{\oplus R}}\Lambda(\D y\times \D
s)\leq \int_{\R_+}(l+R)^2(2lR+\pi R^2)\pr_L(\D l) <\infty
\qquad\forall R>0,$$ so condition \eqref{cond Bool 2} (and
similarly the assumption (A2)) is satisfied. Hence
Theorem~\ref{teo dens} applies and we get
\begin{eqnarray*}\lambda_{\Theta}(x_1,x_2)&\stackrel{\eqref{lambda}}{=}&\int_0^\infty\frac{1}{2\pi}\int_0^{2\pi}
\int_0^l f(x_1-\tau\cos\alpha,x_2-\tau\sin\alpha)\D
\tau\,\D\alpha\,\pr_L(\D l)\\
&=& \int_0^\infty\frac{1}{2\pi}\int_0^{2\pi} \Big((x_1^2+x_2^2)l-
(x_1\cos\alpha+x_2\sin\alpha)l^2+\frac{1}{3}l^3\Big)\D\alpha\,\pr_L(\D
l)\\ &=&(x_1^2+x_2^2)\E[L]+\frac{1}{3}\E[L^3].
\end{eqnarray*}

Note that in the particular case $f\equiv c>0$, $\Theta$ is
stationary and by Corollary~\ref{cor stat} we obtain the well
known result $\lambda_\Theta(x)=c\E[L]$ $\forall x\in\R^d$
(cf.~\cite{Benes-Rataj},\,p.\,42). }\end{example}

\section{On the estimation of the mean density}\label{sec estimation}

While the problem of the mean density estimation has been widely
examined in the stationary case,  the inhomogeneous one has been
mainly faced by assuming local stationarity. For instance,
$\Theta$ is assumed to have a gradient structure, i.e.~it is
considered to be homogeneous perpendicularly to a particular
gradient direction (see \cite{HMS}), or it is assumed  to be
homogeneous in certain subregions of $\R^d$, so that the known
results in the homogeneous case can be applied to estimate a
stepwise approximation of the mean density $\lambda_\Theta$. We
also mention that the stationary Boolean model $\Theta$ is often
assumed to have unknown constant intensity $c>0$ and known mark
distribution $Q$, so that only the parameter $c$ has to be
estimated, being $\lambda_\Theta=c\E_Q[\h^n(Z_0)]$ in this case. A
series of results about the estimation of the intensity $c$
 of the underlying
Poisson point process associated to $\Theta$, related to the
estimation of $\lambda_\Theta$, can be found in \cite{Benes-Rataj}
$\mathcal{x}$3.4 (see also \cite{mrkvicka,SKM}).
\\
 Having now an explicit
formula for the mean density of inhomogeneous Boolean models as in
Assumptions,   in all situations in which it is possible to
estimate the intensity $f$ and the mark distribution $Q$ of the
typical grain $Z_0$, an estimation of $\lambda_\Theta$ could be
obtained. Actually, the estimation of $f$ and $Q$ as well as the
computing of the mean density $\lambda_\Theta(x)$ at a given point
$x\in\R^d$
 might be quite hard. In this section we introduce an estimator
  for the mean density
 $\lambda_\Theta(x)$ of Boolean models as in Assumptions in the
 general case of $f$ and $Q$ unknown.

 In
the proof of Theorem~\ref{teo dens} we have shown, in particular,
 that $\lambda_\Theta$ is also given by the limit in \eqref{scambio}
  for $\h^d$-a.e.~$x\in \R^d$. Noticing
that $\pr(x\in\Theta_{\oplus r})=T_\Theta(B_r(x))$, where
$T_\Theta$ is the  capacity (or hitting) functional of $\Theta$
\cite{matheron}, a natural estimator of $\lambda_\Theta(x)$ can be
given in terms of the \emph{empirical capacity functional} of
$\Theta$, without estimating the intensity $f$ and the
distribution of $Z_0$ separately. We recall that the empirical
capacity functional $\widehat{T}^N_\Xi$ based on an i.i.d. random
sample $\Xi_1,\ldots,\Xi_N$ of a random closed set $\Xi$ is
defined as (see, e.g., \cite{feng})
$$\widehat{T}^N_\Xi(K):=\frac{1}{N}\sum_{i=1}^N \1_{\Xi_i\cap K\neq\emptyset}, \qquad\forall
\mbox{ compact }K\subset\R^d,$$  and that the strong law of large
numbers implies that $\widehat{T}^N_\Xi(K)$ converges almost
surely to $T_\Xi(K)$ for any compact subset $K$ of $\R^d$.
\\
Let $\Theta$ be an inhomogeneous Boolean model in $\R^d$ as in
Assumptions and $\Theta_1,\ldots,\Theta_N$ be an i.i.d. random
sample of $\Theta$; for any fixed $r>0$
 we have that
\beq\E\Big[\frac{\widehat
T^N_\Theta(B_r(x))}{b_{d-n}r^{d-n}}\Big]=\frac{\pr(x\in\Theta_{\oplus
r})}{b_{d-n}r^{d-n}}\ ,\quad{\rm var}\Big(\frac{\widehat
T^N_\Theta(B_r(x))}{b_{d-n}r^{d-n}}\Big)=\frac{\pr(x\in\Theta_{\oplus
r})(1-\pr(x\in\Theta_{\oplus
r}))}{N(b_{d-n}r^{d-n})^2}.\label{attesa stima}\eeq Hence,
\eqref{scambio} and \eqref{attesa stima} suggest that we take
\beq\widehat \lambda_{\Theta}^N(x) :=
\frac{\sum_{i=1}^N\1_{\Theta_i\cap B_{R_N}(x)\neq\emptyset}}{N
b_{d-n}R_N^{d-n}},\label{stimatore}\eeq with $R_N$ such that
\beq\lim_{N\to\infty}R_N= 0 \quad\mbox{and}\qquad
\lim_{N\to\infty}NR_N^{d-n}=\infty,\label{RN}\eeq as
asymptotically unbiased and consistent estimator of the mean
density $\lambda_{\Theta}(x)$ of $\Theta$ at point $x$.
\bp\label{prop stimatore} Let $\Theta$ be a Boolean model in
$\R^d$ as in Assumptions  and $\{\Theta_i\}_{i\in\N}$ be a
sequence of random closed sets i.i.d. as $\Theta$; then
$$\lim_{N\to\infty}\widehat \lambda_{\Theta}^N(x)=\lambda_{\Theta}(x)\quad \mbox{in probability}, \qquad
\mbox{ for }\h^d\mbox{-a.e.~}x\in\R^d. $$ \ep
\begin{proof}
By the definition of $\widehat\lambda_{\Theta}^N$ it follows
$$\lim_{N\to\infty}\E[\widehat\lambda_{\Theta}^N(x)]=\lim_{N\to\infty}
\frac{\pr(x\in\Theta_{\oplus
R_N})}{b_{d-n}R_N^{d-n}}\stackrel{\eqref{scambio}}{=}\lambda_{\Theta}(x)\qquad
\mbox{for }\h^d\mbox{-a.e. }x\in\R^d.
$$
We check now that the variance of $\widehat\lambda_{\Theta}^N(x)$
goes to 0. Since
\begin{itemize} \item[]
$\displaystyle\lim_{N\to\infty}\pr(x\in\Theta_{\oplus
R_N})=\pr(x\in\Theta)=0$ for $\h^d$-a.e.~$x\in\R^d$ because $n<d$,
\item[]
$\displaystyle\lim_{N\to\infty}\frac{\pr(x\in\Theta_{\oplus
R_N})}{b_{d-n}R_N^{d-n}}=\lambda_{\Theta}(x)\in\R$ for
$\h^d$-a.e.~$x\in\R^d$ by \eqref{scambio},
\end{itemize} and
\begin{itemize} \item[] $\displaystyle\lim_{N\to\infty}\frac{1}{NR_N^{d-n}}=0$
by \eqref{RN},\end{itemize} we have that for
$\h^d$-a.e.~$x\in\R^d$
$$
\lim_{N\to\infty}{\rm
var}(\widehat\lambda_{\Theta}^N(x))=\lim_{N\to\infty}\frac{N\pr(x\in\Theta_{\oplus
R_N})(1-\pr(x\in\Theta_{\oplus R_N}))}{(Nb_{d-n}R_N^{d-n})^2}=0.
$$
 Hence  the thesis follows.
\end{proof}

Then, a problem of statistical interest could be to find the
optimal \emph{width} $R_N$ satisfying condition \eqref{RN} which
minimizes the mean squared error of
$\widehat\lambda_{\Theta}^N(x)$
(i.e.~$\E[(\widehat\lambda_{\Theta}^N(x)-\lambda_\Theta(x))^2] $).
To investigate this problem is not the aim of the present paper
and we leave this as open problem for further developments; we
point out here that $\widehat\lambda_\Theta^N$ can be seen as the
generalization to the case of $n$-dimensional random closed sets,
of the well known estimator of the probability density of a random
point, which is a particular 0-dimensional random closed set.
\br[The special case $n=0$]{\rm Even if the particular case $n=0$
can be handle with much more elementary tools, it is easy to check
that if $\Theta=X$ is a random point in $\R^d$ with probability
density function $f_X$, Eq.~\eqref{scambio} holds with
$\lambda_X=f_X$ (it is sufficient to observe that
$\E[\h^0(X\cap\cdot\,)]=\pr(X\in\cdot\,)$) and the estimator
$\widehat{\lambda}_X^N$ turns to be closely related to the well
known definition of histogram (see \cite{PhD} for details). Let us
consider the case of a random variable $X$  with density $f_X$;
then Proposition~\ref{prop stimatore} applies with $d=1$ and
$n=0$, making explicit the correspondence with the usual
\emph{density estimation by means of histograms} (e.g., see
\cite{pestman} $\mathcal{x}$VII.13). Indeed,  if
$\{X_i\}_{i\in\N}$ is a sequence of i.i.d. random variables with
the same distribution of $X$, we define
$$\widehat f_X(x):=\widehat\lambda_X^N(x) \stackrel{\eqref{stimatore}}{=}
\frac{\sum_{i=1}^N\1_{ B_{R_N}(x)}(X_i)}{Nb_1R_N}=\frac{{\rm
card}\{i:X_i\in I_x\}}{N|I_x|},
$$
where $I_x$ is the interval in $\R$ centered in $x$ with length
$|I_x|=2R_N$ with the usual condition
$$|I_x|\longrightarrow 0 \ \mbox{ and }\ N|I_x|\longrightarrow\infty \qquad \mbox{as }N\to\infty.$$
Therefore, statistical problems and techniques related to the
choice of the ``optimal width'' $R_N$ in the general case $N>0$
could be investigated starting from available results for random
variables. }\er

We conclude  observing how the following corollary of
Theorem~\ref{teo dens}, which could lead to further developments
on this topic and on the estimation of the mean density, is
consistent with a couple of available results in literature.
\begin{corollary}\label{cor lim E}
Any Boolean model $\Theta$ as in Assumptions is absolutely
continuous in mean with mean density
\beq\lambda_{\Theta}(x)=\lim_{r\da 0}\frac{\E[{\rm
card}\{(x_i,s_i)\in\Phi\,:\, (x_i+Z_0(s_i))\cap
B_r(x)\neq\emptyset\}]}{b_{d-n}r^{d-n}} \qquad\mbox{for
}\h^d\mbox{-a.e.~}x\in\R^d.\label{dens numero medio}\eeq
\end{corollary}
\begin{proof}
The assertion  follows directly by \eqref{tesi th}, noticing that
$\Lambda({\cal Z}^{x,r})=\E[{\rm card}\{(x_i,s_i)\in\Phi\,:\,
(x_i+Z_0(s_i))\cap B_r(x)\neq\emptyset\}]$.
\end{proof}

\noindent By Proposition 21 in \cite{AKV} we get that for a
locally finite union $\Theta$ of i.i.d.  random closed sets $E_i$
with Hausdorff dimension $n<d$ it holds \beq\lim_{r\da
0}\frac{\pr(x\in\Theta_{\oplus r})}{b_{d-n}r^{d-n}}= \lim_{r\da
0}\frac{\E[{\rm card}\{E_i\,:\,E_i\cap
B_r(x)\neq\emptyset\}]}{b_{d-n}r^{d-n}} \qquad \mbox{for
}\h^d\mbox{-a.e.~}x\in\R^d,\label{lim attesa}\eeq provided that at
least one of the two limits exists. On the other hand, in
\cite{schneider} the mean density of a class of nonstationary
\emph{$n$-flat processes} is studied and a similar result to
\eqref{dens numero medio} is obtained. Namely, we remind that a
$n$-flat process in $\R^d$ (with $n$ integer less than $d$) is a
point process $X$ on ${\cal E}_n^d$, the space of $n$-dimensional
planes in $\R^d$; it is proved that if the intensity measure of
$X$ has a continuous density $h$ with respect to some
translation-invariant, locally finite measure on ${\cal E}_n^d$,
then the $n$-dimensional random closed set $\Theta:=\bigcup_{E\in
X}E$ has continuous mean density
$$\lambda_\Theta(x)=\int_{{\cal L}_n^d} h(x+L)\Psi(\D L),$$
where ${\cal L}_n^d$ is the Grassmannian of $n$-dimensional linear
subspaces in $\R^d$ and $\Psi$ is a finite measure on ${\cal
L}_n^d$ coming from a decomposition result of the intensity
measure of $X$, and  in particular it is claimed that (see
\cite{schneider}, p.\,142, or  \cite{CIME}, p.\,179)
$$\lambda_\Theta(x)=\lim_{r\da 0}\frac{\E[{\rm card}\{E\in
X\,:\,E\cap B_r(x)\neq\emptyset\}]}{b_{d-n}r^{d-n}}.
$$
We may like to notice that  Theorem~\ref{teorema riassuntivo}
applies to $n$-flat process; hence  we are lead to conjecture that
the exchange between limit and integral in \eqref{main2} may hold
for further processes $\Theta=\bigcup_iE_i$, union of i.i.d.
$n$-dimensional random closed sets, so that Proposition~\ref{prop
stimatore} and, by \eqref{lim attesa},  Corollary~\ref{cor lim E},
could be extended to this kinds of random closed sets.
\\

 \noindent{\bf Acknowledgements} Fruitful discussions on the
 Minkowski content notion are acknowledged to Prof.~Luigi
 Ambrosio (SNS, Pisa).

 \end{document}